\documentclass[11pt]{amsart}

\usepackage[latin1]{inputenc}
\usepackage{amsmath,amssymb,latexsym,graphics,color,graphicx,pb-diagram,pb-xy, amscd}
\usepackage[all]{xy}

\usepackage[left=2.54cm, right=2.54cm, top=2.54cm, bottom=2.54cm]{geometry}

\newtheorem{theorem}{Theorem}[section]
\newtheorem{prop}[theorem]{Proposition}
\newtheorem{lem}[theorem]{Lemma}

\newtheorem{definition}[theorem]{Definition}

\newtheorem{remark}[theorem]{Remark}

\newcommand{\af}{\mathrm{af}}
\newcommand{\al}{\alpha}
\newcommand{\alv}[1]{{\alpha^\vee_{#1}}}

\newcommand{\av}{{\alpha^\vee}}

\newcommand{\cA}{\mathcal{A}}

\newcommand{\cl}{\mathrm{cl}}
\newcommand{\clproj}{\cl}

\newcommand{\gaf}{\geh_\af}
\newcommand{\geh}{\mathfrak{g}}
\newcommand{\Hom}{\mathrm{Hom}}

\newcommand{\I}{I}
\newcommand{\Iaf}{I_\af}
\newcommand{\IP}{J}

\newcommand{\la}{\lambda}
\newcommand{\La}{\Lambda}
\newcommand{\lev}{\mathrm{level}}

\newcommand{\om}{\omega}
\newcommand{\OmegaP}{\Omega_J}

\newcommand{\pair}[2]{\langle #1\,,\,#2\rangle}
\newcommand{\Phiaf}{\Phi^\af}
\newcommand{\Phiafm}{\Phi^{\af-}}
\newcommand{\Phiafp}{\Phi^{\af+}}
\newcommand{\PhiafpP}{\Phi^{\af+}_J}

\newcommand{\PhiP}{\Phi_J}

\newcommand{\piP}[1]{\lfloor #1 \rfloor}
\newcommand{\ppiP}{\pi_J}

\newcommand{\qleft}[1]{\xleftarrow{#1}}

\newcommand{\QB}{\mathrm{QB}}
\newcommand{\Qv}{Q^\vee}
\newcommand{\QvP}{Q_J^\vee}
\newcommand{\rhoP}{\rho_J}

\newcommand{\tal}{{\widetilde{\alpha}}}
\newcommand{\ti}[1]{\widetilde{#1}}

\newcommand{\W}{W}
\newcommand{\Waf}{W_\af}
\newcommand{\Wafm}{W_\af^-}

\newcommand{\WP}{W_J}
\newcommand{\WPaf}{(W_J)_\af}
\newcommand{\WUP}{W^J}
\newcommand{\WUPaf}{(W^J)_\af}
\newcommand{\X}{X}
\newcommand{\Xaf}{X_\af}
\newcommand{\Xafv}{X_\af^\vee}
\newcommand{\Xafz}{X_\af^0}

\newcommand{\Z}{\mathbb{Z}}

\newcommand{\BZ}{\mathbb{Z}}
\newcommand{\BQ}{\mathbb{Q}}
\newcommand{\BR}{\mathbb{R}}
\newcommand{\BB}{\mathbb{B}}

\newcommand{\QLS}{\text{\rm QLS}}
\newcommand{\Deg}{\text{\rm Deg}}
\newcommand{\wt}{\text{\rm wt}}

\newcommand{\bi}{{\bf i}}
\newcommand{\bd}{{\bf d}}
\newcommand{\bzero}{{\bf 0}}

\newcommand{\Bpair}[2]{\bigl\langle #1,\,#2 \bigr\rangle}

\newcommand{\ud}[1]{\underline{#1}}

\begin{document}

\title[A uniform model for Kirillov--Reshetikhin crystals]{A uniform model for Kirillov--Reshetikhin crystals\\ {\rm{\small{Extended Abstract}}}}

\author[C.~Lenart]{Cristian Lenart}
\address{Department of Mathematics and Statistics, State University of New York at Albany, 
Albany, NY 12222, U.S.A.}
\email{clenart@albany.edu}
\urladdr{http://www.albany.edu/\~{}lenart/}
%\thanks{Partially supported by the Research in Pairs program of Mathematisches Forschungsinstitut 
%Oberwolfach and the NSF grant DMS--1101264.}

\author[S.~Naito]{Satoshi Naito}
\address{Department of Mathematics, Tokyo Institute of Technology,
2-12-1 Oh-Okayama, Meguro-ku, Tokyo 152-8551, Japan}
\email{naito@math.titech.ac.jp}
%\thanks{Partially supported by the Research in Pairs program of Mathematisches Forschungsinstitut 
%Oberwolfach and Grant-in-Aid for Scientific Research (C), No. 24540010}

\author[D.~Sagaki]{Daisuke Sagaki}
\address{Institute of Mathematics, University of Tsukuba, 
Tsukuba, Ibaraki 305-8571, Japan}
\email{sagaki@math.tsukuba.ac.jp}

\author[A.~Schilling]{Anne Schilling}
\address{Department of Mathematics, University of California, One Shields
Avenue, Davis, CA 95616-8633, U.S.A.}
\email{anne@math.ucdavis.edu}
\urladdr{http://www.math.ucdavis.edu/\~{}anne}
%\thanks{Partially supported by the Research in Pairs program of Mathematisches Forschungsinstitut 
%Oberwolfach and the NSF grants DMS--1001256 and OCI--1147247.}

\author[M.~Shimozono]{Mark Shimozono}
\address{Department of Mathematics, Virginia Polytechnic Institute
and State University, Blacksburg, VA 24061-0123, U.S.A.}
\email{mshimo@vt.edu}


% No need to include the dates
%\received{(Some date $\leq$ 1$^{\mbox{\footnotesize{st}}}$ December 2009)}
%\revised{\today}
%\accepted{tomorrow}

\bibliographystyle{abbrvnat}

\begin{abstract}
We present a uniform
construction of tensor products of one-column Kirillov--Reshetikhin (KR) crystals in
all untwisted affine types, which uses
a generalization of the Lakshmibai--Seshadri paths (in the theory of the 
Littelmann path model). This generalization is based on the graph on parabolic cosets
of a Weyl group known as the parabolic quantum Bruhat graph. A related model is the
so-called quantum alcove model. The proof is based on two lifts of the parabolic quantum Bruhat graph:
to the Bruhat order on the affine 
Weyl group and to Littelmann's poset on level-zero weights. Our construction 
leads to a simple calculation of the energy function. It also implies the 
equality between a Macdonald polynomial specialized at $t=0$ and the graded character of a tensor product of KR modules.
\end{abstract}

\maketitle

%%%%%%%%%%%%%%%%%%%%%%%%%%%%%%%%%%%%%%%%%%%%%%%%%%% 
\section{Introduction}

Our goal in this series of papers (see~\cite{LNSSS.I,LNSSS.II}) is to obtain a uniform construction of tensor 
products of one-column
Kirillov--Reshetikhin (KR) crystals. As a consequence we shall prove
the equality $P_\lambda(q)=X_\lambda(q)$,
where $P_\lambda(q)$ is the Macdonald polynomial $P_\lambda(q,t)$ specialized at $t=0$
and $X_\lambda(q)$ is the graded character of a simple Lie algebra
coming from tensor products of KR modules.
Both the Macdonald polynomials and KR modules are of arbitrary untwisted affine type.
The index $\lambda$ is a dominant weight for the simple Lie subalgebra obtained
by removing the affine node.
Macdonald polynomials and characters of KR modules have been
studied extensively in connection with various fields such as 
statistical mechanics and integrable systems, 
representation theory of Coxeter groups and Lie algebras (and their quantized analogues
given by Hecke algebras and quantized universal enveloping algebras),
geometry of singularities of Schubert varieties, and combinatorics.

Our point of departure is a theorem of Ion~\cite{Ion}, which asserts 
that the nonsymmetric Macdonald polynomials 
at $t=0$ are characters of Demazure submodules of highest weight
modules over affine algebras. This holds for the Langlands duals of 
untwisted affine root systems (and type $A_{2n}^{(2)}$ in the case of 
nonsymmetric Koornwinder polynomials).
Our results apply to the untwisted affine root systems.
The overlapping cases are the simply-laced affine root systems $A_n^{(1)}$,
$D_n^{(1)}$ and $E_{6,7,8}^{(1)}$.

It is known~\cite{FL1,FSS,KMOU,KMOTU,ST,Na} 
that certain affine Demazure characters (including those for the simply-laced affine root systems)
can be expressed in terms of KR crystals, which motivates the relation
between $P$ and $X$. For types $A_n^{(1)}$ and $C_n^{(1)}$, the equality $P=X$
was achieved in~\cite{Le,LeS} by establishing a combinatorial formula for the Macdonald polynomials
at $t=0$ from the Ram--Yip formula~\cite{RY}, and by using explicit models for the one-column KR crystals~\cite{FOS}.
It should be noted that, in types $A_n^{(1)}$ and $C_n^{(1)}$, the one-column KR modules are irreducible
when restricted to the canonical simple Lie subalgebra, while in general this is not the case.
For the cases considered by Ion~\cite{Ion}, the corresponding KR crystals are perfect. This is not necessarily true for the untwisted affine root systems considered in this work, especially for the untwisted non-simply-laced affine root systems.

In this work we provide a type-free approach to the equality $P=X$ for untwisted
affine root systems. 
Lenart's specialization~\cite{Le} of the Ram--Yip formula for Macdonald polynomials uses
the quantum alcove model \cite{LeL}, whose objects are paths in the 
quantum Bruhat graph (QBG), which was defined and studied in \cite{BFP} in relation to the quantum cohomology of the flag variety. On the other hand, Naito and Sagaki~\cite{NS1,NS2,NS4,NS3}
gave models for tensor products of KR crystals of one-column type in terms of projections of level-zero Lakshmibai--Seshadri (LS) paths to the classical weight lattice. Hence we need to establish 
a bijection between the quantum alcove model and projected level-zero LS paths.

In analogy with \cite{BFP} and inspired by the quantum Schubert calculus of homogeneous spaces \cite{Mi,P}
we define the parabolic quantum Bruhat graph (PQBG), 
which is a directed graph structure on parabolic quotients of the Weyl group with respect to a parabolic subgroup.
We construct two lifts of the PQBG. 
The {\em first lift} is from the PQBG to the Bruhat 
order of the affine Weyl group. This is a parabolic analogue of the lift
of the QBG to the affine Bruhat order \cite{LS}, which is the
combinatorial structure underlying Peterson's
theorem \cite{P}; the latter equates the Gromov-Witten invariants of finite-dimensional homogeneous spaces  
with the Pontryagin homology structure constants of Schubert varieties in the affine Grassmannian.
We obtain Diamond Lemmas for the PQBG via projection of the
standard Diamond Lemmas for the affine Weyl group.
We find a {\em second lift} of the PQBG into a poset of Littelmann~\cite{Li} for level-zero weights
and characterize its local structure (such as cover relations) in terms of the PQBG.
Littelmann's poset was defined in connection with LS paths for arbitrary (not necessarily dominant) 
weights, but the local structure was not previously known. The weight poset precisely controls the combinatorics of the
level-zero LS paths and therefore their classical projections, which we formulate directly as quantum LS paths.
Finally, we describe a bijection between the quantum alcove model and the quantum LS paths.

The paper is organized as follows. In Section~\ref{s2} we prepare the background and define the PQBG.
Section~\ref{s3} is reserved for the two lifts of the PQBG and the statement of the Diamond Lemmas.
In Section~\ref{s4} we describe KR crystals in terms of the quantum LS paths and the quantum alcove model in \cite{LeL}; we also give simple combinatorial formulas for the energy function. Finally, we conclude in Section~\ref{s5} with the results on
(nonsymmetric) Macdonald polynomials at $t=0$.

%%%%%%%%%%%%%%%%%%%%%%%%%%%%%%%%%%%%%%%%%%%%%%%%%%%%%%%
\subsection*{Acknowledgments}

The first two and last two authors would like to thank the Mathematisches Forschungsinstitut Oberwolfach for 
their support during the Research in Pairs program, where some of the main ideas of this paper were conceived. 
We would also like to thank Thomas Lam for helpful discussions during FPSAC12 in Nagoya, Japan
and Daniel Orr for his discussions about Ion's work~\cite{Ion}.
We used {\sc Sage}~\cite{sage} and {\sc Sage-combinat}~\cite{sagecombinat} to discover properties about the
level-zero weight poset and to obtain some of the pictures in this paper.

C.L. was partially supported by the NSF grant DMS--1101264.
S.N. was supported by Grant-in-Aid for Scientific Research (C), No. 24540010, Japan.
D.S. was supported by Grant-in-Aid for Young Scientists (B) No.23740003, Japan. 
A.S. was partially supported by the NSF grants DMS--1001256 and OCI--1147247.
M.S. was partially supported by the NSF grant DMS-1200804.

%%%%%%%%%%%%%%%%%%%%%%%%%%%%%%%%%%%%%%%%%%%%%%%%%%%%%%%%
\section{Background}\label{s2}
%%%%%%%%%%%%%%%%%%%%%%%%%%%%%%%%%%%%%%%%%%%%%%%%%%%%%%%%

%%%%%%%%%%%%%%%%%%%%%%%%%%%%%%%%%%%%%%%%%%%%%%%%%%%%%%%%
\subsection{Untwisted affine root datum}

Let $\Iaf=I\sqcup\{0\}$ (resp. $I$) be the Dynkin node set of an untwisted affine algebra $\gaf$ (resp. its canonical subalgebra $\geh$),
$\Waf$ (resp. $\W$) the affine (resp. finite) Weyl group with simple reflections $r_i$ for $i\in \Iaf$
(resp. $i\in I$), and $\Xaf = \Z\delta \oplus\bigoplus_{i\in \Iaf} \Z\La_i$ (resp. $\X = \bigoplus_{i\in I} \Z \om_i$) the affine (resp. finite) weight lattice.
Let $\{\alpha_i \mid i \in \Iaf\}$ be the simple roots, 
$\Phiaf = \Waf \,\{\alpha_i\mid i\in \Iaf\}$ (resp. $\Phi=\W\, \{\alpha_i\mid i\in I\}$) the set of 
affine real roots (resp. roots), and $\Phiafp = \Phiaf \cap \bigoplus_{i\in \Iaf} \Z_{\ge0} \alpha_i$
(resp. $\Phi^+=\Phi\cap \bigoplus_{i\in I} \Z_{\ge0} \alpha_i$) the set of positive affine real (resp. positive) roots.
Furthermore, $\Phiafm=-\Phiafp$ (resp.  $\Phi^-=-\Phi^+$) are the negative affine real (resp. negative) roots.
Let $\Xafv=\Hom_{\Z}(\Xaf,\Z)$ be the dual lattice, $\pair{\cdot}{\cdot}: \Xafv \times \Xaf\to \Z$ the evaluation pairing,
and $\{d\}\cup \{\alv{i} \mid i\in \Iaf \}$ the dual basis of $\Xafv$. 
The natural projection $\cl:\Xaf\to \X$
has kernel $\Z\La_0\oplus\Z\delta$ and sends $\La_i\mapsto\omega_i$ for $i\in I$.

The affine Weyl group $\Waf$ acts on $\Xaf$ and $\Xafv$ by
\begin{equation*}
	r_i \la = \la - \pair{\alpha_i^\vee}{\la} \alpha_i \qquad \text{and} \qquad r_i \mu = \mu - \pair{\mu}{\alpha_i} \alv{i}\,,
\end{equation*}
for $i\in \Iaf$, $\la\in \Xaf$, and $\mu\in \Xafv$.
For $\beta\in \Phiaf$, let $w\in \Waf$ and $i\in\Iaf$ be such that $\beta = w \alpha_i$.
Define the associated reflection $r_\beta\in \Waf$ and associated coroot $\beta^\vee\in \Xafv$ by
$r_\beta = w r_i w^{-1}$ and $\beta^\vee = w \alv{i}$.

The null root is the unique element
$\delta\in \bigoplus_{i\in\Iaf}\Z_{>0} \alpha_i$ which generates the rank 1 sublattice
$\{\la\in \Xaf\mid \pair{\alv{i}}{\la}=0 \text{ for all $i\in \Iaf$} \}$. 
We have $\delta = \alpha_0 + \theta$, where $\theta$ is the highest root for $\geh$.
The canonical central element is the unique
element $c\in \bigoplus_{i\in\Iaf} \Z_{>0} \alv{i}$ which generates the rank 1 sublattice
$\{\mu\in\Xafv\mid \pair{\mu}{\alpha_i}=0 \text{ for all $i\in \Iaf$} \}$.
The level of a weight $\la\in \Xaf$ is defined by $\lev(\la) = \pair{c}{\la}$. 
Let $\Xafz\subset\Xaf$ be the sublattice of level-zero elements.

We denote by $\ell(w)$ for $w\in \Waf$ (resp. $W$) the length of $w$ and by $\lessdot$ the Bruhat cover.
The element $t_\mu\in \Waf$ is the translation by the element $\mu$ in the coroot lattice $Q^\vee=\bigoplus_{i\in I} \Z \alpha_i^\vee$.

%%%%%%%%%%%%%%%%%%%%%%%%%%%%%%%%%%%%%%%%%%%%%%%%%%%%%%%%
\subsection{Affinization of a weight stabilizer}\label{s:stab}

Let $\lambda\in \X$ be a dominant weight, which is fixed throughout the remainder of Sections~\ref{s2} 
and~\ref{s3}.
Let $\WP$ be the stabilizer of $\lambda$ in $\W$.
It is a parabolic subgroup, being generated by $r_i$ for $i\in \IP$, where $\IP = \{i\in \I\mid \pair{\alv{i}}{\lambda}=0\}$.
Let $\QvP = \bigoplus_{i\in \IP} \Z\alv{i}$ be the associated
coroot lattice, 
$\WUP$ the set of minimum-length coset representatives in $\W/\WP$,
$\PhiP\supset \PhiP^+$ the set of roots and positive roots respectively, 
and $\rhoP = (1/2)\sum_{\alpha\in\PhiP^+} \alpha$ (if $J=\emptyset$, then $\rhoP$ is denoted by $\rho$).
Define
\begin{align}
\WPaf&=\WP \ltimes \QvP=\{wt_\mu\in\Waf\mid w\in\WP\,,\mu\in\QvP\}\,,\\
\label{E:Paffinepos}
  \PhiafpP &= \{ \beta\in \Phiafp \mid \cl(\beta) \in \PhiP \}= \PhiP^+ \cup (\Z_{>0}\delta + \PhiP)\,,\\
\label{E:PaffineWpos}
  \WUPaf &= \{x\in \Waf\mid  x\cdot \beta > 0 \text{ for all $\beta\in \PhiafpP$}\, \}.
\end{align}

By \cite[Lemma 10.5]{LS} \cite{P}, any $w\in \Waf$ factors uniquely as $w = w_1 w_2$, where  
$w_1\in \WUPaf$ and $w_2\in \WPaf$. 
Therefore, we can define $\ppiP:\Waf\to \WUPaf$ by $w\mapsto w_1$. We say that $\mu\in \Qv$ is {\em $J$-adjusted} if $\ppiP(t_\mu) = z_\mu t_\mu$ with $z_\mu\in W$. % (in fact, $z_\mu\in W^J$).
Say that $\mu\in \Qv$ is {\em $J$-superantidominant} if $\mu$ is antidominant (i.e., $\pair{\mu}{\alpha} \le 0$ for all $\alpha\in \Phi^+$), and
$\pair{\mu}{\alpha} \ll 0$ for $\alpha\in\Phi^+\setminus\PhiP^+$.

%%%%%%%%%%%%%%%%%%%%%%%%%%%%%%%%%%%%%%%%%%%%%%%%%%%%%%%%%
\subsection{The parabolic quantum Bruhat graph}
The \emph{parabolic quantum Bruhat graph} $\QB(W^J)$ is a directed graph with vertex set $\WUP$,
whose directed edges have the form $w\overset{\al}\longrightarrow \piP{wr_\alpha}$ for $w\in \WUP$ and $\alpha\in \Phi^+\setminus \PhiP^+$.
Here we denote by $\piP{v}$ the minimum-length representative in the coset $v \WP$ for $v\in W$.
There are two kinds of edges:
\begin{enumerate}
\item (Bruhat edge) $w\lessdot wr_\al$. (One may deduce that $wr_\al\in \WUP$.)
\item (Quantum edge) $\ell(\piP{wr_\al})=\ell(w)+1-\pair{\av}{2\rho-2\rhoP}$.
\end{enumerate}
If $J=\emptyset$, then we recover the quantum Bruhat graph $\QB(W)$ defined in \cite{BFP}.

%%%%%%%%%%%%%%%%%%%%%%%%%%%%%%%%%%%%%%%%%%%%%%%%%%%%%%%%
\section{Two lifts of the parabolic quantum Bruhat graph}\label{s3}
%%%%%%%%%%%%%%%%%%%%%%%%%%%%%%%%%%%%%%%%%%%%%%%%%%%%%%%%

%%%%%%%%%%%%%%%%%%%%%%%%%%%%%%%%%%%%%%%%%%%%%%%%%%%%%%%%
\subsection{Lifting $\QB(W^J)$ to $\Waf$}
We construct a parabolic analogue of the lift of $\QB(W)$ to $\Waf$ given in \cite{LS}.  

Let $\OmegaP^\infty\subset\Waf$ be the subset of 
elements of the form $w\ppiP(t_\mu)$
with $w\in \WUP$ and $\mu\in \Qv$  $J$-superantidominant 
 and $J$-adjusted. 
We have $\OmegaP^\infty \subset \WUPaf\cap \Wafm$, where $\Wafm$ is the set of minimum-length coset representatives in $\Waf/\W$. 
Impose the Bruhat covers in $\OmegaP^\infty$ whenever the connecting root 
has classical part in $\Phi\setminus\PhiP$. Then $\OmegaP^\infty$
is a subposet of the Bruhat poset $\Waf$.

\begin{prop} \label{P:Plift} Every edge in $\QB(\WUP)$
lifts to a downward Bruhat cover in $\OmegaP^\infty$, and every
cover in $\OmegaP^\infty$ projects to an edge in $\QB(\WUP)$. More precisely:
\begin{enumerate}
\item
For any edge $\piP{wr_\al}\qleft{\al} w$ in $\QB(\WUP)$,
 and $\mu\in \Qv$ that is $J$-superantidominant
and $J$-adjusted with $\pi_J(t_\mu)=z t_\mu$,
there is a covering relation $y\lessdot x$ in $\OmegaP^\infty$
where 
\[x=wzt_\mu\,,\;\;\;\; y=xr_\tal=wr_\al t_{\chi\av} z t_\mu\,,\;\;\;\;\tal = z^{-1}\al + (\chi + \pair{\mu}{z^{-1}\al})\delta\in \Phiafm\,,\]
and $\chi$ is $0$ or $1$ according as the
arrow in $\QB(\WUP)$ is of Bruhat or quantum type respectively.
\item
Suppose $y\lessdot x$ is an arbitrary covering relation in $\OmegaP^\infty$.
Then we can write $x=wzt_\mu$ with $w\in \WUP$, $z=z_\mu\in \WP$, and $\mu\in \Qv$ $J$-superantidominant and 
$J$-adjusted, as well as $y=xr_\gamma$ with $\gamma=z^{-1}\al+n\delta\in\Phiaf$, $\al\in\Phi^+\setminus\PhiP^+$, and $n\in\Z$.
With the notation $\chi:=n-\pair{\mu}{z^{-1}\al}$, we have 
\[\chi\in\{0,1\}\,,\;\;\;\;\gamma = z^{-1}\al + (\chi+\pair{\mu}{z^{-1}\al})\delta\in \Phiafm\,;\]
furthermore, there is an edge $wr_\al z\qleft{z^{-1}\al} wz$ in $\QB(\W)$ 
and an edge $\piP{wr_\al} \qleft{\al} w$ in $\QB(\WUP)$, where both edges are of Bruhat
type if $\chi=0$ and of quantum type if $\chi=1$.
\end{enumerate}
\end{prop}

%%%%%%%%%%%%%%%%%%%%%%%%%%%%%%%%%%%%%%%%%%%%%%%%%%%%%%%%%
\subsection{The Diamond Lemmas}

In the following, a dotted (resp. plain) edge represents a quantum (resp. Bruhat) edge
in $\QB(\WUP)$, whereas a dashed edge can be of both types. Given $w\in W^J$ and 
$\gamma\in\Phi^+$, define $z\in W_J$ by $r_\theta w = \lfloor r_\theta w \rfloor z$. 
We now state the Diamond Lemmas for $\QB(W^J)$. They are proved based on the lift of $\QB(W^J)$ to $\Waf$ in Proposition \ref{P:Plift} and the fact that such a lemma holds for any Coxeter group~\cite{BB}.

\begin{lem} \label{lemma.diamond.PQBG}
Let $\alpha\in\Phi$ be a simple root, $\gamma \in \Phi^+\setminus\Phi_J^+$, 
and $w\in W^J$. Then we have the following cases, in each of which the bottom two edges imply the top two edges, and vice versa. 
\begin{enumerate}
\item \label{case.1}
In the left diagram we assume $\gamma\ne w^{-1}(\alpha)$. In both cases we have
$r_\alpha \lfloor w r_\gamma\rfloor= \lfloor r_\alpha  w r_\gamma\rfloor$.
\begin{equation}\label{dq1}
	\begin{diagram}\node[2]{r_\alpha \lfloor w r_\gamma\rfloor}\\
	\node{r_\alpha w} \arrow{ne,t}{\gamma} 
	\node[2]{\lfloor w r_\gamma\rfloor} \arrow{nw,t}{\lfloor w r_\gamma\rfloor^{-1}(\alpha)} \\\node[2]{w} \arrow{nw,b}{w^{-1}(\alpha)} \arrow{ne,b}{\gamma}\end{diagram}\qquad\qquad
	\begin{diagram}\node[2]{r_\alpha \lfloor w r_\gamma \rfloor}\\
	\node{r_\alpha w} \arrow{ne,t,..}{\gamma}
	\node[2]{\lfloor w r_\gamma\rfloor} \arrow{nw,t}{\lfloor w r_\gamma \rfloor^{-1}(\alpha)} \\\node[2]{w} \arrow{nw,b}{w^{-1}(\alpha)} \arrow{ne,b,..}{\gamma}\end{diagram}
\end{equation}
\item \label{case.3} 
Here $z$ is defined as above. We assume $\gamma\ne -w^{-1}(\theta)$ whenever both of the hypothesized edges are quantum ones. In the left diagram, the dashed edge is a quantum (resp. a  Bruhat) edge depending on $\langle w^{-1}(\theta),\gamma^\vee\rangle$ being nonzero (resp. zero). In the right diagram, the dashed edge is a Bruhat (resp. a  quantum) edge depending on $\langle w^{-1}(\theta),\gamma^\vee\rangle$ being nonzero (resp. zero).
\begin{equation}\label{dq3}
	\begin{diagram}\node[2]{\lfloor r_\theta w r_\gamma \rfloor}\\
	\node{\lfloor r_\theta w \rfloor} \arrow{ne,t,--}{z(\gamma)}
	\node[2]{\lfloor w r_\gamma \rfloor} \arrow{nw,t,..}{- \lfloor w r_\gamma \rfloor^{-1}(\theta)} \\ \node[2]{w} \arrow{nw,b,..}{-w^{-1}(\theta)} \arrow{ne,b}{\gamma}\end{diagram}
\qquad \qquad
	\begin{diagram}\node[2]{\lfloor r_\theta w r_\gamma \rfloor}\\
	\node{\lfloor r_\theta w \rfloor} \arrow{ne,t,--}{z(\gamma)}
	\node[2]{\lfloor w r_\gamma\rfloor} \arrow{nw,t,..}{-\lfloor w r_\gamma \rfloor^{-1}(\theta)} \\\node[2]{w} \arrow{nw,b,..}{-w^{-1}(\theta)} \arrow{ne,b,..}{\gamma}
	\end{diagram}
\end{equation}
\end{enumerate}
\end{lem}

%%%%%%%%%%%%%%%%%%%%%%%%%%%%%%%%%%%%%%%%%%%%%%%%%%%%%%%%
\subsection{Lifting $\QB(W^J)$ to the level-zero weight poset}
\label{section.level 0 WP}

In~\cite{Li}, Littelmann introduced a poset related to LS paths
for arbitrary (not necessarily dominant) integral weights. 
Littelmann did not give a precise local description of it.
We consider this poset for level-zero weights and
characterize its cover relations in terms of the PQBG.

Let $\lambda\in X$ be a fixed dominant weight (cf. Section \ref{s:stab} and the notation thereof, e.g., $W_J$ is the stabilizer of $\lambda$).  We view $X$ as a sublattice of $\Xafz$. Let $\Xafz(\lambda)$ be the orbit $\Waf \lambda$. 

\begin{definition} {\rm (Level-zero weight poset~\cite{Li})}
\label{definition.level-0 weight poset}
A poset structure is defined on $\Xafz(\lambda)$ as the transitive closure of the relation
\[
 \mu < r_\beta(\mu) \quad \Leftrightarrow \quad	\langle \mu, \beta^\vee \rangle >0\,,
\]
where $\beta\in\Phiafp$. This poset is called the {\em level-zero weight poset for} $\lambda$.
\end{definition}

The cover $\mu\lessdot\nu=r_\beta(\mu)$ of $\Xafz(\lambda)$ is labeled by the root
$\beta\in\Phiafp$. The projection map $\cl$ restricts to the map
$\clproj: \Xafz(\lambda)\to W \lambda$. We identify $W\lambda\simeq W/W_J\simeq W^J$, and consider $\QB(W^J)$.
Our main result is the construction of a lift of $\QB(W^J)$ to $\Xafz(\lambda)$.
The proof is based on Lemma \ref{lemma.diamond.PQBG}.

\begin{theorem} \label{theorem.level 0 weight poset}
Let $\mu\in \Xafz(\lambda)$ and $w:=\clproj(\mu)\in W^J$.
If $\mu \lessdot \nu$ is a cover in $\Xafz(\lambda)$, then its label $\beta$ is in $\Phi^+$ or $\delta-\Phi^+$. Moreover, $w \to \clproj(\nu)$ 
is an up (respectively down) edge in $\QB(W^J)$ labeled by 
$w^{-1}(\beta)\in\Phi^+\setminus\Phi_J^+$ (respectively $w^{-1}(\beta-\delta)$), depending on 
$\beta\in\Phi^+$ (respectively $\beta\in\delta-\Phi^+$). Conversely, if $\begin{diagram}\dgARROWLENGTH=1.5em\node{w} \arrow{e,t}{\gamma}\node{wr_\gamma=w'}\end{diagram}$ 
(respectively $\begin{diagram}\dgARROWLENGTH=1.5em\node{w} \arrow{e,t,..}{\gamma}
\node{\lfloor wr_\gamma\rfloor=w'}\end{diagram}$) in $\QB(W^J)$ for
$\gamma \in \Phi^{+} \setminus \Phi_{J}^{+}$, then there exists a cover $\mu \lessdot \nu$ in 
$\Xafz(\lambda)$ labeled by $w(\gamma)$ (respectively $\delta+w(\gamma)$) with $\clproj(\nu)=w'$. 
\end{theorem}

%%%%%%%%%%%%%%%%%%%%%%%%%%%%%%%%%%%%%%%%%%%%%%%%
\section{Models for KR crystals and the energy function} \label{s4}
%%%%%%%%%%%%%%%%%%%%%%%%%%%%%%%%%%%%%%%%%%%%%%%%

%%%%%%%%%%%%%%%%%%%%%%%%%%%%%%
\subsection{Quantum LS paths}
%%%%%%%%%%%%%%%%%%%%%%%%%%%%%%

Throughout this section, 
we fix a dominant integral weight $\lambda \in X$, and set 
$J:=\bigl\{i \in I \mid \pair{\alpha_{i}^{\vee}}{\lambda}=0\bigr\}$. 

\begin{definition} \label{dfn:QBG-achain}
Let $x,\,y \in W^{J}$, and let $\sigma \in \BQ$ 
be such that $0 < \sigma < 1$. A directed $\sigma$-path 
from $y$ to $x$ is, by definition, a directed path 
\begin{equation*}
x=w_{0} \stackrel{\gamma_{1}}{\leftarrow} w_{1}
\stackrel{\gamma_{2}}{\leftarrow} w_{2} 
\stackrel{\gamma_{3}}{\leftarrow} \cdots 
\stackrel{\gamma_{n}}{\leftarrow} w_{n}=y
\end{equation*}
from $y$ to $x$ in $\QB(W^J)$ such that
$\sigma \pair{\gamma_{k}^{\vee}}{\lambda} \in \BZ$ 
for all $1 \le k \le n$.
\end{definition}

A quantum LS path of shape $\lambda$ is 
a pair $\eta=(\ud{x}\,;\,\ud{a})$ of 
a sequence $\ud{x}\,:\,x_{1},\,x_{2},\,\dots,\,x_{s}$ of 
elements in $W^{J}$ with $x_{u} \ne x_{u+1}$ 
for $1 \le u \le s-1$ and a sequence 
$\ud{\sigma}\,:\,
 0=\sigma_{0} < \sigma_{1} < \cdots < \sigma_{s}=1$ 
of rational numbers such that 
there exists a directed $\sigma_{u}$-path 
from $x_{u+1}$ to $x_{u}$ for each $1 \le u \le s-1$. 
Denote by $\QLS(\lambda)$ the set of quantum LS paths 
of shape $\lambda$. 
We identify an element 
$\eta=(x_{1},\,x_{2},\,\dots,\,x_{s}\,;\,
\sigma_{0},\,\sigma_{1},\,\dots,\,\sigma_{s}) \in \QLS(\lambda)$
with the following piecewise linear, continuous map 
$\eta:[0,1] \rightarrow \BR \otimes_{\BZ} X$:
\begin{equation*}
\eta(t)=\sum_{u'=1}^{u-1}
(\sigma_{u'}-\sigma_{u'-1})x_{u'} \cdot \lambda + 
(t-\sigma_{u-1})x_{u} \cdot \lambda \quad \text{for \ }
\sigma_{u-1} \le t \le \sigma_{u}, \  1 \le u \le s,
\end{equation*}
and set $\wt(\eta)=:\eta(1)$. Following \cite{Li}, we define the root operators 
$e_{i}$ and $f_{i}$ for $i \in I_{\af}=I \sqcup \{0\}$ 
as follows. For $\eta \in \QLS(\lambda)$ and 
$i \in I_{\af}$, we set
\begin{equation*}
\begin{array}{l}
H(t)=H^{\eta}_{i}(t):=\pair{\alpha_{i}^{\vee}}{\eta(t)}
  \quad \text{for \,} t \in [0,1], \\[3mm]
m=m^{\eta}_{i}
  :=\min\bigl\{H^{\eta}_{i}(t) \mid t \in [0,1]\bigr\};
\end{array}
\end{equation*}
in fact, $m \in \BZ_{\le 0}$. 
If $m = 0$, then $e_{i}\eta:=\bzero$. 
If $m \le -1$, 
then we define $e_{i}\eta$ by:
\begin{equation*}
(e_{i}\eta)(t)=
\begin{cases}
\eta(t) & \text{if } 0 \le t \le t_{0}, \\[2mm]
r_{i,\,m+1}(\eta(t))=
\eta(t_{0})+s_{i}(\eta(t)-\eta(t_{0}))
       & \text{if \,} t_{0} \le t \le t_{1}, \\[2mm]
r_{i,\,m+1}r_{i,\,m}(\eta(t))=
\eta(t)+\ti{\alpha}_{i} & \text{if \,} t_{1} \le t \le 1,
\end{cases}
\end{equation*}
where
\begin{equation*}
\begin{array}{l}
t_{1}:=\min\bigl\{t \in [0,1] \mid 
       H(t)=m \bigr\}, \\[2mm]
t_{0}:=\max\bigl\{t \in [0,t_{1}] \mid
       H(t) = m+1 \bigr\},
\end{array}
\end{equation*}
$r_{i,\,n}$ is the reflection with respect to the hyperplane 
$H_{i,\,n}:=\bigl\{\mu \in \BR \otimes_{\BZ} X \mid 
 \pair{\alpha_{i}^{\vee}}{\mu}=n \bigr\}$ 
for each $n \in \BZ$, and 
\begin{equation*}
\ti{\alpha}_{i}:=
\begin{cases}
\alpha_{i} & \text{if $i \ne 0$}, \\[1.5mm]
-\theta & \text{if $i = 0$},
\end{cases}
\qquad
s_{i}:=
\begin{cases}
r_{i} & \text{if $i \ne 0$}, \\[1.5mm]
r_{\theta} & \text{if $i = 0$}.
\end{cases}
\end{equation*}

\begin{figure}[htbp]
\begin{center}
  \scalebox{.8}{\includegraphics{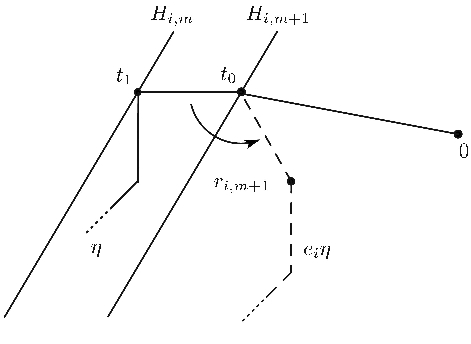}}
  \caption{Root operator $e_{i}$.}
\end{center}
\end{figure}
The definition of $f_{i}\eta$ is similar. 
The following theorem is one of our main results.

%
%%%%%%%%%%%%%%%%
%%% thm:LScl %%%
%%%%%%%%%%%%%%%%
%
\begin{theorem} \label{thm:LScl}
\begin{enumerate}
\item[(1)]
The set $\QLS(\lambda)$ together with crystal operators $e_{i}$, $f_{i}$ for
$i \in I_{\af}$ and weight function $\wt$, becomes a regular crystal 
with weight lattice $X$.
\item[(2)]
For each $i \in I$, the crystal $\QLS(\om_{i})$ 
is isomorphic to the crystal basis of $W(\om_{i})$, 
the fundamental representation of level zero, 
introduced by Kashiwara {\rm \cite{Kas-OnL}}.
\item[(3)]
Let $\bi=(i_{1},\,i_{2},\,\dots,\,i_{p})$ be an arbitrary sequence of 
elements of $I$, and set 
$\lambda_{\bi}:=\om_{i_{1}}+\om_{i_{2}}+ \cdots + \om_{i_{p}}$. 
There exists a crystal isomorphism 
$\Psi_{\bi}:\QLS(\lambda_{\bi}) \stackrel{\sim}{\rightarrow} 
\QLS(\om_{i_{1}}) \otimes \QLS(\om_{i_{2}}) \otimes 
\cdots \otimes \QLS(\om_{i_{p}})$.
\end{enumerate}
\end{theorem}
%
%%%%%%%%%%%%%%%%
%%% rem:LScl %%%
%%%%%%%%%%%%%%%%
%

\begin{remark} \label{rem:LScl}{\rm 
It is known that 
the fundamental representation $W(\om_{i})$ (of level zero) 
is isomorphic to the KR module $W_{1}^{(i)}$ 
in the sense of {\rm \cite[\S2.3]{HKOTT}}
(but for the explicit form of the Drinfeld polynomials
of $W(\om_{i})$, see {\rm \cite[Remark~3.3]{N}}), 
and that it has a global crystal basis (see {\rm \cite[Theorem 5.17]{Kas-OnL}}).
Furthermore, the crystal basis of 
$W(\om_{i}) \cong W_{1}^{(i)}$ is unique, 
up to a nonzero constant multiple 
(see also {\rm \cite[{Lemma} 1.5.3]{NS-CMP}});
we call it a (one-column) KR crystal. 
By the theorem above, the crystal $\QLS(\lambda)$ is a model 
for the corresponding tensor product of KR crystals.
}
\end{remark}

%%%%%%%%%%%%%%%%%%%%%%%%%%%%%%%%%%%%%%%%%%%%%%%%%%%%%%%
\subsection{Sketch of the proof of Theorem~{\rm \ref{thm:LScl}}}
%%%%%%%%%%%%%%%%%%%%%%%%%%%%%%%%%%%%%%%%%%%%%%%%%%%%%%%

First, let us recall the definition of LS paths 
of shape $\lambda$ from \cite{Li}. 

\begin{definition}
For $\mu,\,\nu \in X_{\af}^{0}(\lambda)$ with $\mu > \nu$ 
(see Definition {\rm \ref{definition.level-0 weight poset}}) and a rational number $0 < \sigma < 1$, 
a $\sigma$-chain for $(\mu,\nu)$ is, 
by definition, a sequence $\mu=\xi_{0} \gtrdot \xi_{1} \gtrdot
\dots \gtrdot \xi_{n}=\nu$ of covers in $X_{\af}^{0}(\lambda)$ 
such that $\sigma\pair{\gamma_{k}^{\vee}}{\xi_{k-1}} \in \BZ$ 
for all $k=1,\,2,\,\dots,\,n$, 
where $\gamma_{k}$ is the label for 
$\xi_{k-1} \gtrdot \xi_{k}$. 
\end{definition}

\begin{definition}
An LS path of shape $\lambda$ is, by definition, 
a pair $(\ud{\nu}\,;\,\ud{\sigma})$ of a sequence 
$\ud{\nu}:\nu_{1} > \nu_{2} > \cdots > \nu_{s}$ of 
elements in $X_{\af}^{0}(\lambda)$ and a sequence 
$\ud{\sigma}:0=\sigma_{0} < \sigma_{1} < \cdots < \sigma_{s}=1$ of 
rational numbers such that 
there exists a $\sigma_{u}$-chain for $(\nu_{u},\,\nu_{u+1})$ 
for each $u=1,\,2,\,\dots,\,s-1$. 
\end{definition}

We denote by $\BB(\lambda)$ the set of all LS paths of shape $\lambda$.
We identify an element 
$$\pi=(\nu_{1},\,\nu_{2},\,\dots,\,\nu_{s}\,;\,
\sigma_{0},\,\sigma_{1},\,\dots,\,\sigma_{s}) \in \BB(\lambda)$$
with the following piecewise linear, continuous map 
$\eta:[0,1] \rightarrow \BR \otimes_{\BZ} X_{\af}$:
\begin{equation*}
\pi(t)=\sum_{u'=1}^{u-1}
(\sigma_{u'}-\sigma_{u'-1})\nu_{u'}+
(t-\sigma_{u-1})\nu_{u} \quad \text{for \ }
\sigma_{u-1} \le t \le \sigma_{u}, \  1 \le u \le s.
\end{equation*}
Let
$\BB(\lambda)_{\cl}:=\bigl\{\cl(\pi) \mid \pi \in \BB(\lambda) \bigr\}$,
where $\cl(\pi)$ is the piecewise linear, continuous map 
$[0,1]\to \BR \otimes_{\BZ} X$ defined by
$(\cl(\pi))(t):=\cl(\pi(t))$ for $t \in [0,1]$.
For $\eta \in \BB(\lambda)_{\cl}$ and $i \in I_{\af}$, 
we define $e_{i}\eta$ and $f_{i}\eta$ in exactly 
the same way as above. Then it is known from 
\cite{NS1, NS2} that the same statement as in 
Theorem~\ref{thm:LScl} holds for $\BB(\lambda)_{\cl}$. 
Thus, Theorem~\ref{thm:LScl} follows immediately from 
the following proposition. 

\begin{prop}
The affine crystals $\QLS(\lambda)$ and $\BB(\lambda)_{\cl}$ are isomorphic.
\end{prop}

This proposition is a consequence of Theorem~\ref{theorem.level 0 weight poset}. 
Let us show that if $\eta \in \QLS(\lambda)$, then 
$\eta \in \BB(\lambda)_{\cl}$. Here, for simplicity, 
we assume that $\eta=(x,\,y\,;\,0,\,
\sigma,\,1) \in \QLS(\lambda)$, and 
$x=w_{0} \stackrel{\gamma_{1}}{\leftarrow} w_{1}
\stackrel{\gamma_{2}}{\leftarrow} w_{2}=y$ is a 
directed $\sigma$-path from $y$ to $x$.
Take an arbitrary $\mu \in X_{\af}^{0}(\lambda)$ 
such that $\cl(\mu)=y$. By applying Theorem~\ref{theorem.level 0 weight poset}
to $w_{1} \stackrel{\gamma_{2}}{\leftarrow} w_{2}=y$, 
we obtain a cover $\nu_{1} \gtrdot \mu$ 
for some $\nu_{1} \in X_{\af}^{0}(\lambda)$ with 
$\cl(\nu_{1})=w_{1}$. Then, by applying Theorem~\ref{theorem.level 0 weight poset}
to $x=w_{0} \stackrel{\gamma_{1}}{\leftarrow} w_{1}$, 
we obtain a cover $\nu \gtrdot \nu_{1}$ 
for some $\nu \in X_{\af}^{0}(\lambda)$ with 
$\cl(\nu)=x$. Thus we get a sequence 
$\nu \gtrdot \nu_{1} \gtrdot \mu$ of 
covers in $X_{\af}^{0}(\lambda)$. It can be easily seen that 
this is a $\sigma$-chain for $(\nu,\,\mu)$, which implies that 
$\pi=(\nu,\,\mu \,;\, 0,\,\sigma,\,1) \in \BB(\lambda)$. 
Therefore, $\eta=\cl(\pi) \in \BB(\lambda)_{\cl}$. 
The reverse inclusion can be shown similarly. 

%%%%%%%%%%%%%%%%%%%%%%%%%%%%%%%%%%%%%%%%%%%%%%%%
\subsection{Description of the energy function in terms of quantum LS paths}
%%%%%%%%%%%%%%%%%%%%%%%%%%%%%%%%%%%%%%%%%%%%%%%%

Recall the notation in Theorem~\ref{thm:LScl}\,(2); 
for simplicity, we set $\lambda:=\lambda_{\bi}$. 
In \cite{NS3}, Naito and Sagaki introduced 
a degree function $\Deg_{\lambda}:\BB(\lambda)_{\cl}=\QLS(\lambda)
\rightarrow \BZ_{\le 0}$, and proved that $\Deg_{\lambda}$ is 
identical to the {\em energy function} \cite{HKOTT,HKOTY} on the tensor product 
$\BB(\om_{i_{1}})_{\cl} \otimes \BB(\om_{i_{2}})_{\cl} \otimes 
\cdots \otimes \BB(\om_{i_{p}})_{\cl}=
\QLS(\om_{i_{1}}) \otimes \QLS(\om_{i_{2}}) \otimes 
\cdots \otimes \QLS(\om_{i_{p}})$ (which is isomorphic to 
the corresponding tensor product of KR crystals; see Remark~\ref{rem:LScl}) 
via the isomorphism $\Psi_{\bi}$. 
The function $\Deg_{\lambda}:\BB(\lambda)_{\cl} 
\rightarrow \BZ_{\le 0}$ is described in terms of $\QB(W^J)$ as follows.
For $x,\,y \in W^{J}$ let
\begin{equation*}
\bd : 
x=
 w_{0} \stackrel{\beta_{1}}{\leftarrow}
 w_{1} \stackrel{\beta_{2}}{\leftarrow} \cdots 
       \stackrel{\beta_{n}}{\leftarrow}
 w_{n}=y
\end{equation*}
be a shortest directed path from $y$ to $x$, and define
\begin{equation*}
\wt(\bd):=\sum_{
 \begin{subarray}{c}
 1 \le k \le n \text{ such that } \\[1mm]
 \text{$w_{k-1} \stackrel{\beta_{k}}{\leftarrow} w_{k}$
 is a down arrow}
 \end{subarray}
}
\beta_{k}^{\vee}. 
\end{equation*}
The value $\pair{\wt(\bd)}{\lambda}$ does not 
depend on the choice of a shortest directed path $\bd$ 
from $y$ to $x$, and 

\begin{theorem}
Let $\eta=(x_{1},\,x_{2},\,\dots,\,x_{s}\,;\,
\sigma_{0},\,\sigma_{1},\,\dots,\,\sigma_{s}) \in 
\QLS(\lambda)=\BB(\lambda)_{\cl}$. Then, 
\begin{equation}\label{eq:deg}
\Deg(\eta)=-\sum_{u=1}^{s-1} 
 (1-\sigma_{u})\Bpair{\lambda}{\wt(\bd_{u})},
\end{equation}
where $\bd_{u}$ is a shortest directed path from $x_{u+1}$ to $x_{u}$.
\end{theorem}

%%%%%%%%%%%%%%%%%%%%%%%%%%%%%%%%%%%%%%%%%%%%%%%%%%%%%%%%
\subsection{The quantum alcove model}\label{S:qalc} 
The quantum alcove model is a generalization of the alcove model of the first author and Postnikov \cite{LP,LP1}, which, 
in turn, is a discrete counterpart of the Littelmann path model \cite{Li}. For the affine Weyl group terminology below, we refer to \cite{H}.
Fix a dominant weight $\lambda$. We say that two alcoves are {adjacent} if they are distinct and have a 
common wall. Given a pair of adjacent alcoves $A$ and $B$, we write 
$A \stackrel{\beta}{\longrightarrow} B$  for $\beta\in\Phi$ if the common wall is orthogonal to $\beta$ and 
$\beta$ points in the direction from $A$ to $B$.

\begin{definition}\label{def:lch}{\rm \cite{LP}}
		The sequence of roots $(\beta_1, \beta_2, \dots, \beta_m)$ is called a
		\emph{$\lambda$-chain} if 
		\[	
		A_0=A_{\circ} \stackrel{-\beta_1}{\longrightarrow} A_1
			%\stackrel{\beta_2}{\longrightarrow} A_2
		\stackrel{-\beta_2}{\longrightarrow}\dots 
		\stackrel{-\beta_m}{\longrightarrow} A_m=A_\circ -\lambda\]
is a shortest sequence of alcoves from the fundamental alcove $A_\circ$ to its translation by $-\lambda$.
\end{definition}

The {\em lex $\lambda$-chain} $\Gamma_{\mathrm{lex}}$ is a particular $\lambda$-chain
defined in \cite[Section 4]{LP1}. Given an arbitrary $\lambda$-chain $\Gamma=(\beta_1, \beta_2, \dots, \beta_m)$, 
let $r_i=r_{\beta_i}$ and define the {\em level sequence} $(l_1,\ldots,l_m)$ of $\Gamma$ by
$l_i=|\left\{ j \geq i \, \mid \, \beta_j = \beta_i \right\} |$.

\begin{definition}\label{def:admissible}{\rm \cite{LeL}}
A (possibly empty) finite subset $J=\left\{ j_1 < j_2 < \cdots < j_s \right\}$ of $\{1,\ldots,m\}$
 is a \emph{$\Gamma$-admissible subset} if
we have the following path in $\QB(W)$:
\begin{equation}
	\label{eqn:admissible}
	1 \stackrel{\beta_{j_1}}{\longrightarrow} r_{j_1} \stackrel{\beta_{j_2}}{\longrightarrow} r_{j_1}r_{j_2} 
	\stackrel{\beta_{j_3}}{\longrightarrow} \cdots \stackrel{\beta_{j_s}}{\longrightarrow} r_{j_1}r_{j_2}\cdots r_{j_s}=\kappa(J)\,.
\end{equation}
Let $\cA^\Gamma(\lambda)$ be the collection of $\Gamma$-admissible subsets.
Define the level of $J\in \cA^\Gamma(\lambda)$ by 
${\rm level}(J)= \sum_{j \in J^{-}} {l}_j$, where $J^-\subseteq J$ corresponds to the down steps in Bruhat order in {\rm \eqref{eqn:admissible}}. 
\end{definition}

\begin{theorem} \label{T:qalcove2qls} There is an isomorphism of graded classical crystals
$\Xi:\cA^{\Gamma_{\mathrm{lex}}}(\lambda)\to \QLS(\lambda)$.
\end{theorem}
The map $\Xi$ is the following forgetful map.
Given the path \eqref{eqn:admissible}, based on the structure of $\Gamma_{\rm lex}$
we select a subpath, and project its elements under $W\to W/W_\lambda$ where $W_\lambda$ is the
stabilizer of $\lambda$ in $W$, thereby obtaining a quantum LS path. The inverse map is more subtle, and is 
based on the so-called {\em tilted Bruhat theorem} of \cite{LNSSS.I}; this is a QBG analogue of the 
minimum-length Deodhar lift~\cite{D} $W/W_\lambda\to W$.

An affine crystal structure was defined on $\cA^{\Gamma_{\rm lex}}(\lambda)$ in \cite{LeL}. We show that $\Xi$
is an affine crystal isomorphism, up to removing some $f_0$-arrows from $\QLS(\lambda)$. The remaining arrows are called {\em Demazure arrows} in \cite{ST}, as they correspond to 
the arrows of a certain affine Demazure crystal, cf. \cite{FSS}. Moreover the above bijection
sends the level statistic to the $\Deg$ statistic of \eqref{eq:deg}.
Thus, we have proved that $\cA^{\Gamma_{\rm lex}}(\lambda)$ is also a model for KR crystals. 

In \cite{LeL2}, we show that all the affine crystals $\cA^\Gamma(\lambda)$, for various $\Gamma$, 
are isomorphic. Making particular choices in classical types, we can translate the level statistic into a 
so-called {\em charge statistic} on sequences of the corresponding {\em Kashiwara--Nakashima columns} \cite{KN}. 
In type $A$, we recover the classical Lascoux--Sch\"utzenberger charge \cite{LSc}. Type $C$ was 
worked out in \cite{Le,LeS}, while type $B$ is considered in \cite{BL}. 

%%%%%%%%%%%%%%%%%%%%%%%%%%%%%%%%%%%%%%%%%%%%%%%%%%%%%%%%%%
\section{Macdonald polynomials} \label{s5}
%%%%%%%%%%%%%%%%%%%%%%%%%%%%%%%%%%%%%%%%%%%%%%%%%%%%%%%%%%

The symmetric Macdonald polynomials $P_\lambda(x;q,t)$ \cite{macaha} are a remarkable family of orthogonal 
polynomials associated to any affine root system (where $\lambda$ is a dominant weight for the canonical
finite root system), which depend on 
parameters $q,t$. They generalize the corresponding irreducible characters, which are recovered upon setting $q=t=0$.
For untwisted affine root systems the Ram--Yip formula \cite{RY} expresses $P_\lambda(x;q,t)$ in terms of all subsequences of any $\lambda$-chain $\Gamma$
(cf. Definition \ref{def:lch}). In \cite{Le}, it was shown that the Ram--Yip formula takes the following simple form 
for $t=0$:
\begin{equation}\label{eq:ry}
P_{\lambda}(x;q,0)=\sum_{J\in\cA^\Gamma(\lambda)} q^{{\rm level}(J)}\,x^{\wt(J)}\,,
\end{equation}
where $\wt(J)$ is a weight associated with $J$. By using the results in 
Section~\ref{S:qalc}, we can write the right-hand side of~\eqref{eq:ry} as a sum over the corresponding tensor 
product of KR crystals. It follows that $P_{\lambda}(x;q,0)=X_\lambda(q)$.
Furthermore, we can express the nonsymmetric Macdonald polynomial $E_{w \lambda}(x;q,0)$, 
for $w \in W$, in a similar way to~\eqref{eq:ry}, by summing over those $J\in\cA^\Gamma(\lambda)$ 
with $\kappa(J)\le w$ (in Bruhat order), where $\kappa(J)$ was defined in~\eqref{eqn:admissible}. This formula 
can be derived both by induction, based on Demazure operators, and from the Ram--Yip formula 
(in this case, for nonsymmetric Macdonald polynomials); however, the latter derivation is more involved 
than the one of~\eqref{eq:ry}, as it uses the transformation on admissible subsets defined in~\cite[Section 5.1]{Le0}.

%%%%%%%%%%%%%%%%%%%%%%%%%%%%%%%%%%%%%%%%%%%%%%%%%%% 

\end{document}